\ifdefined\pdftexversion\pdfoutput=1\fi
\documentclass[11pt]{amsart}

\usepackage[T1]{fontenc}
\ifdefined\pdftexversion\usepackage[utf8]{inputenc}\fi
\usepackage{lmodern}
\usepackage{microtype}

\usepackage[letterpaper, margin=1.2in]{geometry}
\usepackage{amssymb}
\usepackage{mathtools}
\numberwithin{equation}{section}

\usepackage{graphicx}
\usepackage{xcolor}

\definecolor{linknavy}{HTML}{15477E}

\usepackage{hyperref}
\hypersetup{
  colorlinks = true,
  linkcolor  = linknavy,
  citecolor  = linknavy,
  urlcolor   = linknavy,
  pdftitle   = {Chains of Maximum Length in Ornamentation Lattices},
  pdfauthor  = {Minh-Chau Bach, Hyunseo Cho, Xiang Lu, Luka Todorovic},
  pdfsubject = {Combinatorics},
  pdfkeywords= {Ornamentation lattice, spine,
                linear extension, shifted Young tableau},
}

\theoremstyle{definition}
\newtheorem{theorem}{Theorem}[section]
\newtheorem{lemma}[theorem]{Lemma}
\newtheorem{proposition}[theorem]{Proposition}
\newtheorem{corollary}[theorem]{Corollary}
\newtheorem{observation}[theorem]{Observation}
\newtheorem*{thmheight}{Theorem~\ref{thm:height-function}}
\newtheorem*{thmspine}{Theorem~\ref{thm:genetic}}
\newtheorem*{corlinext}{Corollary~\ref{cor:linext}}

\newtheorem{definition}[theorem]{Definition}
\newtheorem{question}[theorem]{Question}

\newcommand{\sfT}{\mathsf{T}}
\newcommand{\sfP}{\mathsf{P}}
\newcommand{\sfS}{\mathsf{S}}
\newcommand{\sfB}{\mathsf{B}}
\newcommand{\sfD}{\mathsf{D}}
\newcommand{\sfY}{\mathsf{Y}}
\newcommand{\sfSp}{\mathsf{Sp}}
\newcommand{\calO}{\mathcal{O}}
\newcommand{\calX}{\mathcal{X}}
\newcommand{\calD}{\mathcal{D}}
\newcommand{\calB}{\mathcal{B}}
\newcommand{\dep}{\operatorname{depth}}
\newcommand{\bd}{\operatorname{bd}}
\newcommand{\desc}[1]{\operatorname{desc}(#1)}
\newcommand{\dn}[1]{{\downarrow}#1}
\newcommand{\ML}{\#\mathrm{ML}}
\newcommand{\LE}{e}
\newcommand{\Trs}{\operatorname{Trs}}
\newcommand{\Spine}{\operatorname{Spine}}

\title{Chains of Maximum Length in Ornamentation Lattices}

\author{Minh-Chau Bach}
\address{Center for Data Science, New York University, New York, NY 10011, USA}
\email{mb8490@nyu.edu}

\author{Hyunseo Cho}
\address{Department of Mathematics, University of Michigan, Ann Arbor, MI 48109, USA}
\email{choamy@umich.edu}

\author{Xiang Lu}
\address{Department of Mathematics, University of Michigan, Ann Arbor, MI 48109, USA}
\email{jameslx@umich.edu}

\author{Luka Todorovic}
\address{Department of Mathematics, University of Michigan, Ann Arbor, MI 48109, USA}
\email{lukat@umich.edu}

\date{November 28, 2026}


\begin{document}

\begin{abstract}
    Ornamentation lattices of finite rooted trees, introduced by Defant and Sack, generalize the Tamari and Boolean lattices.
    We identify the spine of an ornamentation lattice with the lattice of order ideals of a poset of proper ancestor--descendant pairs, thereby reducing the enumeration of chains of maximum length to counting linear extensions.
    We obtain explicit formulas for the number of these chains for trees in which every vertex with at least two children has depth at most~$2$.
\end{abstract}
\maketitle

\section{Introduction and main results}\label{sec:intro}

\begin{figure}[ht]
\centering
\includegraphics{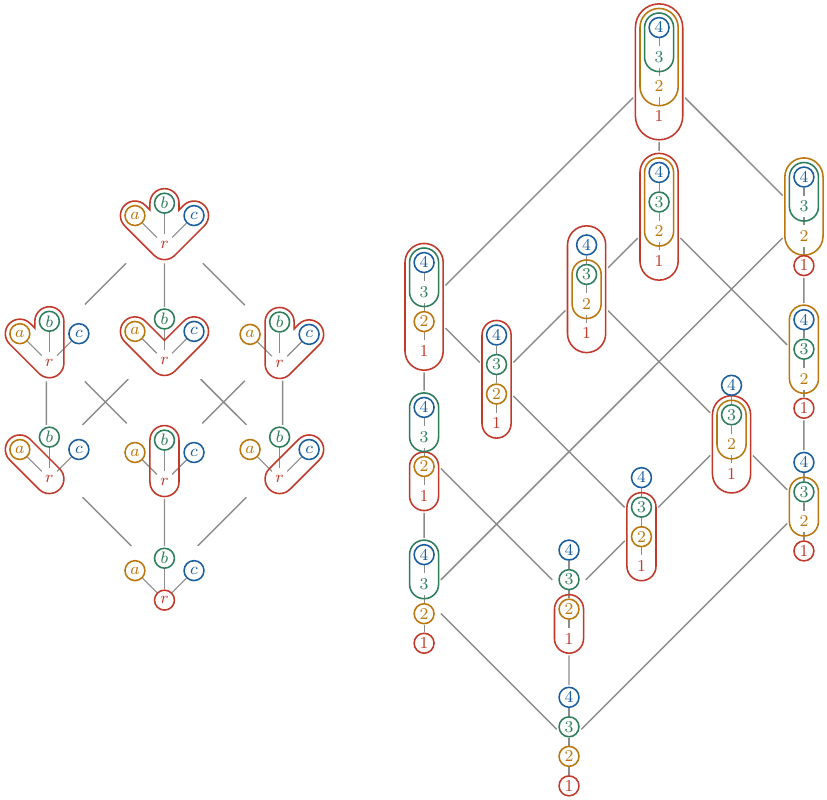}
\caption{The ornamentation lattice of a star with~$3$ leaves (left) is a Boolean lattice, and that of a path with~$4$ vertices (right) is a Tamari lattice.
In each ornamentation, the ornament at a vertex~$v$ is outlined in the color of~$v$.}
\label{fig:specializations}
\end{figure}

Defant and Sack~\cite[Theorem~1.1 and Section~3]{defant2024operahedronlattices} introduced ornamentation lattices to prove that Laplante-Anfossi's posets on the vertices of operahedra~\cite[Section~2.1]{laplante2022diagonal} are lattices.
An \emph{ornamentation}~$\rho$ of a finite rooted tree~$\sfT$ assigns to each vertex~$v$ a subtree~$\rho(v)$ rooted at~$v$, called the ornament at~$v$, such that any two ornaments are nested or disjoint (Definition~\ref{def:ornamentation}).
The ornamentation lattice~$\calO(\sfT)$ is the poset of ornamentations of~$\sfT$ ordered by pointwise containment, with minimum element~$\hat0$ and maximum element~$\hat1$.
For example, the ornamentation lattice of a path rooted at an endpoint is isomorphic to a Tamari lattice~\cite{Tamari62}, and that of a star rooted at its center is isomorphic to a Boolean lattice; see~\cite[Examples~2.6 and~2.7]{abram2025ornamentation} and Figure~\ref{fig:specializations}.
Following a suggestion of Defant and Sack~\cite[Section~8]{defant2024operahedronlattices} to generalize results on Tamari lattices, such as those of Fishel and Nelson~\cite{FishelNelson2014}, we study chains of maximum length in ornamentation lattices of finite rooted trees.

Along a maximal chain from~$\hat0$ to~$\hat1$, the ornament at each vertex~$u$ grows from~$\{u\}$ to the full subtree rooted at~$u$.
As ornaments are nested or disjoint, a cover may enlarge an ornament by more than one vertex.
Consequently, maximal chains can have different lengths (Figure~\ref{fig:lattice-Y}).
However, along a chain of maximum length, each cover adds exactly one vertex to one ornament (Lemma~\ref{lem:one-increment}).

The \emph{spine} is the subposet of~$\calO(\sfT)$ induced by the ornamentations lying on chains of maximum length~\cite[Section~4]{Thomas2006}.
Since~$\calO(\sfT)$ is trim~\cite[Remark~6.3]{ajran2025popstack}, its spine is a distributive sublattice~\cite[Lemma~7]{Thomas2006} (see also Proposition~\ref{prop:lattice-transport}).
To each ornamentation on the spine, we associate the set of pairs~$(u,v)$ for which the ornament at~$u$ contains a proper descendant~$v$ of~$u$.
The resulting sets are precisely the order ideals of the poset~$\calX(\sfT)$ of proper ancestor--descendant pairs, ordered componentwise with ancestors below descendants (Theorem~\ref{thm:genetic}).
Reading the pairs in the order in which they are added along a chain of maximum length yields a linear extension of~$\calX(\sfT)$; conversely, every linear extension arises from exactly one chain of maximum length (Corollary~\ref{cor:linext}).

Fishel and Nelson~\cite[Definition~3.1 and Theorem~3.3]{FishelNelson2014} constructed a bijection between standard~$m$-shifted tableaux and chains of maximum length in the~$m$-Tamari lattices of Bergeron and Pr\'eville-Ratelle~\cite{BergeronPrevilleRatelle12}.
For~$m=1$, the tableaux are standard shifted Young tableaux of staircase shape, so Schur's product formula enumerates the chains~\cite{Schur1911,Thrall52}.
When~$\sfT$ is a path,~$\calX(\sfT)$ is isomorphic to a shifted staircase (Proposition~\ref{prop:paths}).
Corollary~\ref{cor:linext} gives a bijection between chains of maximum length and linear extensions of~$\calX(\sfT)$ for every finite rooted tree~$\sfT$.

\subsection{Structural results}
Throughout the paper,~$\sfT$ is a finite rooted tree with vertex set~$V$ and root~$r$.

\begin{thmheight}
Let~$\rho$ be an ornamentation of~$\sfT$.
The height of~$\rho$ in~$\calO(\sfT)$, that is, the maximum length of a chain from~$\hat0$ to~$\rho$, is
\[
h(\rho)=\sum_{v\in V}\bigl(|\rho(v)|-1\bigr)\,.
\]
\end{thmheight}

Taking~$\rho=\hat1$, we recover Ajran and Defant's result that the maximum chain length equals the number of proper ancestor--descendant pairs in~$\sfT$~\cite[Remark~6.3]{ajran2025popstack}.

Write~$\Spine(\calO(\sfT))$ for the spine of~$\calO(\sfT)$ (Definition~\ref{def:chain-height-spine}), and~$J(P)$ for the lattice of order ideals of a finite poset~$P$.
For an ornamentation~$\rho$, let~$\Phi(\rho)\coloneqq\{(u,v)\in\calX(\sfT)\colon v\in\rho(u)\}$ (Definition~\ref{def:ancestor-descendant}).
By Theorem~\ref{thm:height-function},~$|\Phi(\rho)|$ is the height of~$\rho$.
We call an ornamentation~$\rho$ \emph{genetic} if, for every nonroot vertex~$w$ with parent~$u$, either~$\rho(w)=\{w\}$ or~$\rho(w)\subseteq\rho(u)$ (Definition~\ref{def:genetic}).

\begin{thmspine}
The elements of~$\Spine(\calO(\sfT))$ are exactly the genetic ornamentations.
The map~$\Phi$ restricts to an order isomorphism
\[
\Phi\colon\Spine(\calO(\sfT))\longrightarrow J(\calX(\sfT))\,.
\]
\end{thmspine}

If~$\rho$ is not genetic, then some nonroot vertex~$w$ with parent~$u$ has a nonsingleton ornament~$\rho(w)$ disjoint from~$\rho(u)$.
On every chain from~$\rho$ to~$\hat1$, when the ornaments at~$u$ and~$w$ first intersect, the ornament at~$u$ must gain every vertex already in the ornament at~$w$ (see the proof of Lemma~\ref{lem:maxchain-implies-genetic}).
Thus~$h$ increases by at least~$2$ in one step, so~$\rho$ does not lie on the spine.

\begin{corlinext}
The chains of maximum length in~$\calO(\sfT)$ are in bijection with the linear extensions of~$\calX(\sfT)$.
\end{corlinext}

Ajran and Defant proved that~$\calO(\sfT)$ is semidistributive and extremal~\cite[Theorem~6.2 and Remark~6.3]{ajran2025popstack} (see also~\cite[Theorem~4.12]{abram2025ornamentation}).
Thomas and Williams~\cite[Proposition~2.4]{ThomasWilliams2019} showed that the spine of an extremal lattice is isomorphic to the lattice of order ideals of an associated poset called the Galois poset.
For~$\calO(\sfT)$, Theorem~\ref{thm:genetic} identifies the spine explicitly as~$J(\calX(\sfT))$, so the Galois poset of~$\calO(\sfT)$ is isomorphic to~$\calX(\sfT)$~\cite[Proposition~3.4.2]{StanleyEC1}.

We also realize~$\calO(\sfT)$ as an interval in a lattice of transfer systems (Proposition~\ref{prop:transfer}).
Luo and Rognerud characterized the spine of a lattice of transfer systems~\cite[Proposition~7.3]{LuoRognerud26} and showed that the spine is isomorphic to a lattice of order ideals~\cite[Proposition~8.3]{LuoRognerud26}.
In Section~\ref{sec:transfer}, we use their characterization to give a second proof that every ornamentation on the spine is genetic (Lemma~\ref{lem:maxchain-implies-genetic}).

Other recent work on ornamentation lattices includes that of Abram, Bastidas, G\'elinas, Pilaud, and Sack~\cite[Definition~4.9 and Theorem~4.12]{abram2025ornamentation}, who described the irreducible elements of~$\calO(\sfT)$, and that of Sack~\cite{sack2026generalized}, who generalized ornamentation lattices to pointed building sets.
Beyond ornamentation lattices, Corteel, Jang, and Rognerud~\cite[Theorem~7.12]{CJR26} identified the spine of the extra slow Tamari lattice with the lattice of order ideals of a poset of colored intervals.
In Cambrian lattices~\cite[Section~1.2]{ThomasWilliams2019} and~$(P,\phi)$-Tamari lattices~\cite[Proposition~4.7]{SegoviaTamari25}, chains of maximum length are also in bijection with linear extensions, whereas in lollipop-graph tubing lattices they are in bijection with partially shiftable tableaux~\cite[Theorem~5.18]{DahlbergFishel24}.
\subsection{Enumerative results}
We write~$H(\sfT)$ for the maximum length of a chain in~$\calO(\sfT)$ and~$\ML(\sfT)$ for the number of chains of that length.
By Theorem~\ref{thm:height-function}, we have~$H(\sfT)=\sum_{v\in V}\dep(v)$, since each vertex~$v$ has~$\dep(v)$ proper ancestors.
The \emph{branching depth}~$\bd(\sfT)$ is the greatest depth of a vertex with at least two children, or~$0$ if no such vertex exists (Definition~\ref{def:branching-depth}).
We obtain a product formula for~$\ML(\sfT)$ when~$\bd(\sfT)\le1$ (Corollary~\ref{cor:bd-one}) and a more general formula when~$\bd(\sfT)\le2$ (Theorem~\ref{thm:bd-two}).

As an example, consider the \emph{broom}~$\sfB_{2,k}$, obtained from a two-edge path rooted at one endpoint by attaching~$k$ leaves to the other (Definition~\ref{def:broom}).
By Corollary~\ref{cor:medium-brooms}, for every~$k\ge0$,
\[
\ML(\sfB_{2,k})
=\sum_{q=0}^{k}\frac{k!}{(k-q)!}\,
  \frac{(3k-q)!}{2^q\,6^{k-q}}\,.
\]
A general formula for trees of branching depth at least~$3$ remains open, even for brooms (Question~\ref{q:broom-hook}).

After reviewing the necessary background in Section~\ref{sec:background}, we prove the height formula in Section~\ref{sec:height}.
Section~\ref{sec:maximum-length-genetic} describes the spine and relates~$\calO(\sfT)$ to transfer systems.
Section~\ref{sec:enumeration} decomposes~$\calX(\sfT)$ into blocks.
In Sections~\ref{sec:bd-one} and~\ref{sec:bd-two}, we use this decomposition to enumerate chains of maximum length for trees of branching depth at most~$1$ and~$2$, respectively.
\section{Background}\label{sec:background}
Defant and Sack~\cite[Section~3]{defant2024operahedronlattices} assigned ornaments to the nonroot vertices of a rooted plane tree.
We also assign an ornament to the root, as in~\cite[Definition~2.1]{abram2025ornamentation}.

Define a partial order~$\preceq_\sfT$ on~$V$ by~$u\preceq_\sfT v$ if~$u$ lies on the path from the root~$r$ to~$v$; we write~$u\prec_\sfT v$ when~$u\preceq_\sfT v$ and~$u\ne v$.
If~$u\preceq_\sfT v$, then~$u$ is an \emph{ancestor} of~$v$, and~$v$ is a \emph{descendant} of~$u$.
We use \emph{proper ancestor} and \emph{proper descendant} when~$u\ne v$.
If~$v$ covers~$u$, then~$u$ is the \emph{parent} of~$v$, and~$v$ is a \emph{child} of~$u$.
The \emph{depth} of~$v$, denoted by~$\dep(v)$, is the number of edges on the path from~$r$ to~$v$.
The depth of~$\sfT$ is~$\max_{v\in V}\dep(v)$.
We write~$\desc{u}\coloneqq\{v\in V\colon u\preceq_\sfT v\}$.
For elements~$x\preceq y$ of a poset~$(P,\preceq)$, we write~$[x,y]_P\coloneqq\{z\in P\colon x\preceq z\preceq y\}$.
For the poset~$(V,\preceq_\sfT)$, we use the notation~$[u,v]_\sfT$.
In our drawings, the root is placed at the bottom, with children above their parents.

For~$n\ge1$, let~$\sfP_n$ be the path on vertices~$1,\dots,n$, rooted at~$1$, with edges between consecutive vertices.
For~$n\ge0$, let~$\sfS_n$ be the star whose root has~$n$ children.
We write~$\sfY$ for the rooted tree on~$4$ vertices in which the root~$r$ has a unique child~$a$, which in turn has two children~$b$ and~$c$.

\begin{definition}\label{def:ornament}
    An \emph{ornament} of~$\sfT$ is a nonempty subset of~$V$ that induces a connected subgraph of~$\sfT$.
    An ornament consisting of a single vertex is a \emph{singleton}.
\end{definition}

\begin{definition}\label{def:ornamentation}
    An \emph{ornamentation} of~$\sfT$ is a map~$\rho \colon V \to 2^{V}$ that assigns an ornament to each vertex and satisfies the following two conditions:
    \begin{enumerate}
    \item \emph{Anchoring.} For every~$v \in V$, we have~$v \in \rho(v)$ and~$\rho(v) \subseteq \desc{v}$.
    \item \emph{Laminarity.} For all~$u, v \in V$, the ornaments~$\rho(u)$ and~$\rho(v)$ are nested or disjoint; that is,
          \[
              \rho(u) \subseteq \rho(v) \quad \text{or} \quad
              \rho(v) \subseteq \rho(u) \quad \text{or} \quad
              \rho(u) \cap \rho(v) = \varnothing\,.
          \]
    \end{enumerate}
\end{definition}

By anchoring and the connectedness of ornaments,~$\rho(v)$ induces a subtree rooted at~$v$.
For~$u\prec_\sfT w$, anchoring gives~$u\in\rho(u)$ and~$u\notin\rho(w)$.
Thus laminarity implies
\begin{equation}\label{eq:nested-or-disjoint}
    \rho(w)\subseteq\rho(u)\qquad\text{or}\qquad\rho(u)\cap\rho(w)=\varnothing\,.
\end{equation}

\begin{definition}\label{def:ornamentation-poset}
    Let~$\calO(\sfT)$ be the set of ornamentations of~$\sfT$.
    Define a partial order~$\preceq_{\calO}$ on~$\calO(\sfT)$ by
    \[
    \rho \preceq_{\calO} \sigma
     \quad\text{if}\quad
        \rho(v) \subseteq \sigma(v) \qquad\text{for all } v \in V\,.
    \]
    We write~$\rho\prec_{\calO}\sigma$ for the strict order.
\end{definition}

The poset~$(\calO(\sfT),\preceq_{\calO})$ has minimum~$\hat0$ and maximum~$\hat1$, given by~$\hat0(v)=\{v\}$ and~$\hat1(v)=\desc{v}$ for~$v\in V$.
Defant and Sack~\cite[Proposition~3.1]{defant2024operahedronlattices} proved that this poset is a lattice, the \emph{ornamentation lattice} of~$\sfT$ (see also~\cite[Theorem~2.4]{abram2025ornamentation}).

\begin{definition}\label{def:chain-height-spine}
    A \emph{chain} of a poset is a subset whose elements are pairwise comparable.
    The \emph{length} of a finite chain is its number of elements minus one.
    In a poset with minimum element~$\hat0$, the \emph{height} of an element~$x$ is the maximum length of a chain from~$\hat0$ to~$x$.
    The \emph{spine}~$\Spine(P)$ of a finite poset~$P$ is the subposet induced by the elements that lie on chains of maximum length in~$P$.
\end{definition}

\begin{definition}\label{def:order-ideal}
    A subset~$I\subseteq P$ of a poset~$(P,\preceq)$ is an \emph{order ideal} (or simply an \emph{ideal}) if~$x\in I$ and~$y\preceq x$ imply~$y\in I$.
    Let~$J(P)$ denote the poset of order ideals of~$P$, ordered by inclusion.
\end{definition}

For every finite poset~$P$, the poset~$J(P)$ is a distributive lattice under union and intersection, graded by cardinality~\cite[Section~3.4]{StanleyEC1}.

\section{The height function}\label{sec:height}

We compute the height of an ornamentation by counting the proper descendants in each ornament.

\begin{definition}\label{def:ancestor-descendant}
    Let
    \[
        \calX(\sfT) \coloneqq  \{(u,v)\in V\times V \colon u \prec_\sfT v\}
    \]
    be the set of proper ancestor--descendant pairs of~$\sfT$, ordered componentwise by
    \[
        (u,v)\preceq_{\calX} (u',v')\quad\text{if}\quad u\preceq_\sfT u' \quad\text{and}\quad v\preceq_\sfT v'\,.
    \]
    We write~$(u,v)\prec_{\calX}(u',v')$ for the strict order.
    For a map~$\rho\colon V\to 2^{V}$ with~$v\in\rho(v)\subseteq\desc{v}$ for all~$v\in V$, let
    \[
        \Phi(\rho) \coloneqq  \{(u,v)\in\calX(\sfT)\colon v\in\rho(u)\}\,.
    \]
\end{definition}

Define~$h\colon\calO(\sfT)\to\mathbb{Z}_{\ge0}$ by~$h(\rho)\coloneqq|\Phi(\rho)|$.

Each vertex~$u$ contributes to~$\Phi(\rho)$ the~$|\rho(u)|-1$ pairs~$(u,v)$ with~$v\in\rho(u)\setminus\{u\}$, so
\begin{equation}\label{eq:h-closed-form}
    h(\rho)=\sum_{u\in V}\bigl(|\rho(u)|-1\bigr)\,.
\end{equation}
Since~$\Phi(\hat0)=\varnothing$ and~$\Phi(\hat1)=\calX(\sfT)$, we have~$h(\hat0)=0$ and~$h(\hat1)=|\calX(\sfT)|$.

\begin{lemma}[{cf.~\cite[Section~2]{ajran2025popstack}}]\label{lem:strict-monotonicity}
    The function~$h$ is strictly order-preserving: if~$\sigma \prec_{\calO} \rho$, then~$h(\sigma) < h(\rho)$.
\end{lemma}

\begin{proof}
    We have~$|\sigma(v)|\le|\rho(v)|$ for every~$v\in V$, with strict inequality for at least one vertex, so~$h(\sigma)<h(\rho)$ by~\eqref{eq:h-closed-form}.
\end{proof}

\begin{lemma}[{cf.~\cite[Lemma~2.4]{ajran2025popstack} and~\cite[Lemma~4.3]{abram2025ornamentation}}]\label{lem:predecessor-existence}
    Let~$\rho\in\calO(\sfT)$ satisfy~$h(\rho)>0$.
    Then there exists~$\sigma\in\calO(\sfT)$ with~$\sigma\prec_{\calO}\rho$ and~$h(\sigma)=h(\rho)-1$.
\end{lemma}

\begin{proof}
    Choose a~$\preceq_\sfT$-maximal vertex~$w$ among those with nonsingleton ornaments, and let~${x\ne w}$ be a leaf of the subtree induced by~$\rho(w)$.
    Define~$\sigma\colon V\to2^V$ by~$\sigma(w)\coloneqq\rho(w)\setminus\{x\}$ and~$\sigma(v)\coloneqq\rho(v)$ for~$v\ne w$.
    Deleting the leaf~$x$ preserves connectedness and anchoring.
    For laminarity, let~$z\ne w$.
    If~$\rho(z)$ is disjoint from~$\rho(w)$ or contains~$\rho(w)$, then~$\rho(z)$ is disjoint from~$\sigma(w)$ or contains~$\sigma(w)$.
    Otherwise~$\rho(z)\subseteq\rho(w)$, so~$w\prec_\sfT z$, and the maximality of~$w$ forces~$\rho(z)=\{z\}$, which is contained in or disjoint from~$\sigma(w)$.
    Thus~$\sigma\in\calO(\sfT)$ and~$\sigma\prec_\calO\rho$.
    Only one ornament loses one vertex, so~$h(\sigma)=h(\rho)-1$.
\end{proof}

\begin{theorem}\label{thm:height-function}
    Let~$\rho \in \calO(\sfT)$.
    The height of~$\rho$ in~$\calO(\sfT)$ is~$h(\rho)$.
\end{theorem}

\begin{proof}
    Along every chain~$\hat0=\rho_0\prec_\calO\cdots\prec_\calO\rho_k=\rho$, Lemma~\ref{lem:strict-monotonicity} gives~$k\le h(\rho)$.
    For the reverse inequality, applying Lemma~\ref{lem:predecessor-existence} repeatedly gives a chain of length~$h(\rho)$ from~$\rho$ down to an ornamentation with~$h=0$, which is~$\hat0$.
\end{proof}

\begin{corollary}[{\cite[Remark~6.3]{ajran2025popstack}}]\label{cor:chain-length}
    The maximum length of a chain in~$\calO(\sfT)$ is
    \[
    H(\sfT) = h(\hat1) = |\calX(\sfT)| = \sum_{v \in V} \bigl(\lvert\desc{v}\rvert - 1\bigr)\,.
    \]
\end{corollary}

\begin{proof}
    Apply Theorem~\ref{thm:height-function} to~$\hat1$.
\end{proof}

\begin{figure}[t]
\centering
\includegraphics{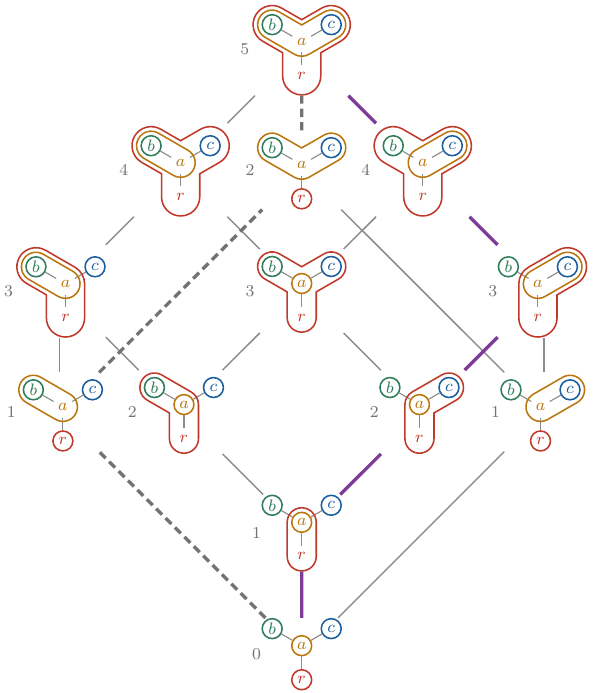}
\caption{Heights in~$\calO(\sfY)$.
The dashed maximal chain has length~$3$, and the purple solid one has length~$5$.}
\label{fig:lattice-Y}
\end{figure}

Figure~\ref{fig:lattice-Y} displays the ornamentation lattice~$\calO(\sfY)$.

\section{Genetic ornamentations and the distributive spine}\label{sec:maximum-length-genetic}
\subsection{The spine}

\begin{definition}\label{def:genetic}
An ornamentation~$\rho$ is \emph{genetic} if, for each vertex~$u$ and each child~$w$ of~$u$, either~$\rho(w)=\{w\}$ or~$\rho(w)\subseteq\rho(u)$.
\end{definition}

\begin{lemma}\label{lem:one-increment}
On a chain of maximum length in~$\calO(\sfT)$, each cover adds one vertex to one ornament.
\end{lemma}

\begin{proof}
Let~$\hat0=\rho_0\prec_{\calO}\rho_1\prec_{\calO}\cdots\prec_{\calO}\rho_N=\hat1$ be the chain.
By Lemma~\ref{lem:strict-monotonicity},
\[
    \varnothing=\Phi(\rho_0)\subsetneq\Phi(\rho_1)\subsetneq\cdots\subsetneq\Phi(\rho_N)=\calX(\sfT)\,,
\]
and~$N=|\calX(\sfT)|$ by Corollary~\ref{cor:chain-length}, so~$|\Phi(\rho_i)\setminus\Phi(\rho_{i-1})|=1$ for~$1\le i\le N$.
\end{proof}

\begin{lemma}\label{lem:maxchain-implies-genetic}
Every element of~$\Spine(\calO(\sfT))$ is genetic.
\end{lemma}

\begin{proof}
Let~$\hat0=\rho_0\prec_{\calO}\cdots\prec_{\calO}\rho_N=\hat1$ be a chain of maximum length, and suppose for contradiction that some~$\rho_i$ is not genetic.
Choose a vertex~$u$ and a child~$w$ of~$u$ such that~$\rho_i(w)$ is a nonsingleton not contained in~$\rho_i(u)$.
By~\eqref{eq:nested-or-disjoint}, the ornaments~$\rho_i(w)$ and~$\rho_i(u)$ are disjoint.

Let~$j>i$ be the first index for which~$\rho_j(w)$ and~$\rho_j(u)$ intersect.
Such an index exists because~$\rho_N(u)=\desc{u}$.
At that step,~\eqref{eq:nested-or-disjoint} gives~$\rho_j(w)\subseteq\rho_j(u)$.
Since~$\rho_{j-1}(u)$ and~$\rho_{j-1}(w)$ are disjoint,
\[
    \rho_i(w)\subseteq\rho_{j-1}(w)
    \subseteq\rho_j(u)\setminus\rho_{j-1}(u)\,.
\]
Thus the ornament at~$u$ gains at least two vertices from~$\rho_{j-1}$ to~$\rho_j$, contradicting Lemma~\ref{lem:one-increment}.
\end{proof}

\begin{definition}\label{def:reconstruction-map}
For~$A\subseteq\calX(\sfT)$, define~$\Psi(A)\colon V\to2^V$ by
\[
    \Psi(A)(u)=\{u\}\cup\{v\in V\colon (u,v)\in A\}
    \qquad(u\in V)\,.
\]
\end{definition}

The map~$\Psi(A)$ satisfies anchoring but need not be an ornamentation.
For every ornamentation~$\rho$ and every~$A\subseteq\calX(\sfT)$, the maps~$\Phi$ and~$\Psi$ satisfy
\begin{equation}\label{eq:recording-reconstruction}
    \Psi(\Phi(\rho))=\rho\,,
    \qquad \Phi(\Psi(A))=A\,.
\end{equation}

\begin{lemma}\label{lem:recording-ideals}
The maps~$\Phi$ and~$\Psi$ restrict to mutually inverse, order-preserving maps between the genetic ornamentations of~$\sfT$ and~$J(\calX(\sfT))$.
\end{lemma}

\begin{proof}
Let~$\rho$ be genetic, and suppose that~$(u',v')\preceq_{\calX}(u,v)\in\Phi(\rho)$.
Since~$\rho(u)$ is a nonsingleton and~$u'\preceq_\sfT u$, iterating the genetic condition gives~$\rho(u)\subseteq\rho(u')$, so~$v\in\rho(u')$.
Connectedness then gives~$v'\in\rho(u')$ since~$u'\prec_\sfT v'\preceq_\sfT v$.
Thus~$(u',v')\in\Phi(\rho)$, proving that~$\Phi(\rho)$ is an ideal.

Conversely, let~$I$ be an order ideal of~$\calX(\sfT)$ and set~$\rho\coloneqq\Psi(I)$.
If~$(u,v)\in I$, then~${(u,z)\in I}$ for every~$z\in[u,v]_\sfT\setminus\{u\}$, so each set~$\rho(u)$ is connected.
If~$u\prec_\sfT w$ and~$\rho(w)$ is a nonsingleton, then for every~$x\in\rho(w)\setminus\{w\}$, both~$(u,w)$ and~$(u,x)$ lie below~$(w,x)\in I$.
Hence~$\rho(w)\subseteq\rho(u)$.
This proves the genetic condition, and also laminarity when~$u$ and~$w$ are comparable, since a singleton~$\rho(w)=\{w\}$ is contained in or disjoint from~$\rho(u)$.
If~$u$ and~$w$ are incomparable, then~$\desc{u}\cap\desc{w}=\varnothing$, so~$\rho(u)$ and~$\rho(w)$ are disjoint.
Thus~$\rho$ is a genetic ornamentation.
Finally, by~\eqref{eq:recording-reconstruction}, the maps are inverse, and both preserve inclusion by definition.
\end{proof}

\begin{figure}[ht]
\centering
\includegraphics{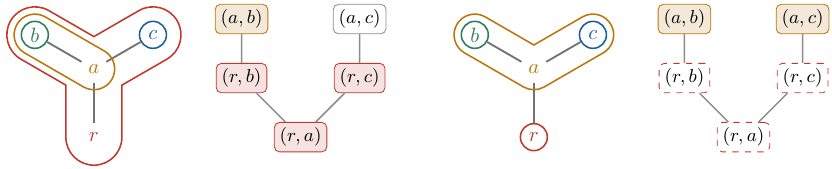}
\caption{A genetic ornamentation (left) and a nongenetic ornamentation (right).
For each ornamentation~$\rho$, filled boxes indicate the pairs in~$\Phi(\rho)$.
Dashed boxes mark the pairs that are below a pair of~$\Phi(\rho)$ but are not in~$\Phi(\rho)$.}
\label{fig:genetic-ideal}
\end{figure}

For the nongenetic ornamentation~$\rho$ of Figure~\ref{fig:genetic-ideal},~$\Phi(\rho)$ contains~$(a,b)$ but omits the smaller pair~$(r,b)$, so it is not an ideal.

\begin{theorem}\label{thm:genetic}
An ornamentation lies in~$\Spine(\calO(\sfT))$ if and only if it is genetic.
The map~$\Phi$ restricts to an order isomorphism
\[
    \Phi\colon\Spine(\calO(\sfT))\longrightarrow J(\calX(\sfT))\,,
    \qquad \Phi^{-1}=\Psi\,.
\]
\end{theorem}

\begin{proof}
By Lemmas~\ref{lem:maxchain-implies-genetic} and~\ref{lem:recording-ideals}, every ornamentation~$\rho$ on the spine is genetic and~$\Phi(\rho)$ is an ideal.
Conversely, let~$\rho$ be genetic, so that~$\Phi(\rho)\in J(\calX(\sfT))$ by Lemma~\ref{lem:recording-ideals}.
Choose a maximal chain~$C$ through~$\Phi(\rho)$ in~$J(\calX(\sfT))$.
Since~$J(\calX(\sfT))$ is graded by cardinality,~$C$ has length~$|\calX(\sfT)|$.
Applying~$\Psi$ to~$C$ gives a chain of the same length through~$\rho$, so~$\rho$ lies on the spine by Corollary~\ref{cor:chain-length}.
The isomorphism now follows from Lemma~\ref{lem:recording-ideals}.
\end{proof}

\begin{figure}[t]
\centering
\includegraphics{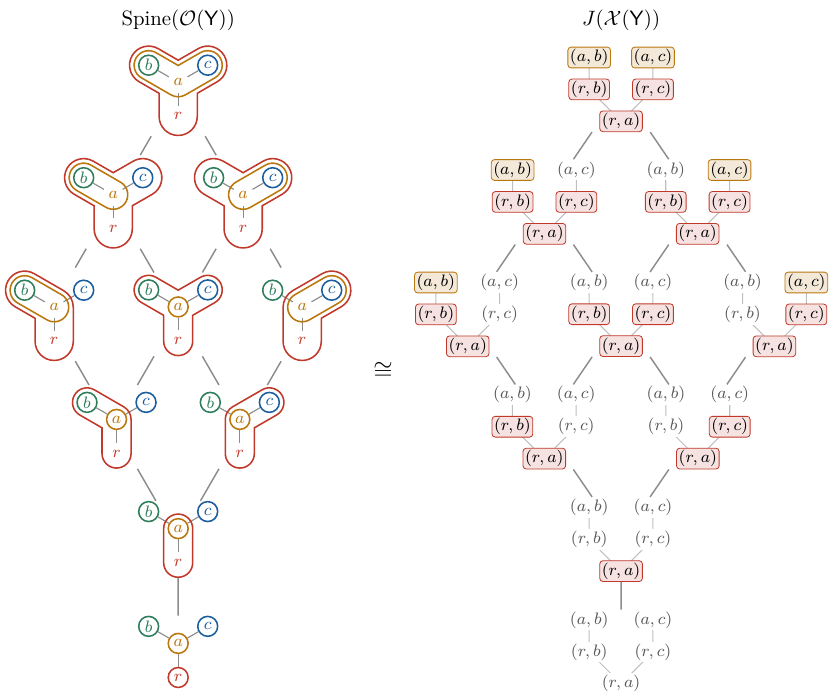}
\caption{The lattices~$\Spine(\calO(\sfY))$ and~$J(\calX(\sfY))$, identified by~$\Phi$.
On the right, colored boxes mark the pairs in each ideal, and a pair~$(u,v)$ is colored like the ornament at~$u$ on the left.}
\label{fig:tree-and-pairposet}
\end{figure}

\begin{corollary}\label{cor:linext}
The chains of maximum length in~$\calO(\sfT)$ are in bijection with the linear extensions of~$\calX(\sfT)$.
\end{corollary}

\begin{proof}
By Theorem~\ref{thm:genetic}, the map~$\Phi$ identifies the chains of maximum length in~$\calO(\sfT)$ with the maximal chains of~$J(\calX(\sfT))$, which correspond to the linear extensions of~$\calX(\sfT)$~\cite[Proposition~3.5.2]{StanleyEC1}.
\end{proof}

Figure~\ref{fig:tree-and-pairposet} shows the~$10$ elements of~$\Spine(\calO(\sfY))$.
The~$6$ maximal chains of~$\Spine(\calO(\sfY))$ are the chains of maximum length in~$\calO(\sfY)$.

\begin{samepage}
\begin{proposition}[{cf.~\cite[Lemma~7]{Thomas2006}}]\label{prop:lattice-transport}
The poset~$\Spine(\calO(\sfT))$ is a distributive sublattice of~$\calO(\sfT)$, graded with rank function~$h$ and containing~$\hat0$ and~$\hat1$.
For~$\rho,\sigma\in\Spine(\calO(\sfT))$ and~$u\in V$, we have
\[
    (\rho\wedge\sigma)(u)=\rho(u)\cap\sigma(u)\,,
    \qquad
    (\rho\vee\sigma)(u)=\rho(u)\cup\sigma(u)\,.
\]
\end{proposition}
\end{samepage}

\begin{proof}
Let~$\rho,\sigma\in\Spine(\calO(\sfT))$.
By Theorem~\ref{thm:genetic},~$\Phi(\rho)$ and~$\Phi(\sigma)$ are ideals, and so are their intersection and union.
Applying~$\Psi$ gives the maps~$u\mapsto\rho(u)\cap\sigma(u)$ and~$u\mapsto\rho(u)\cup\sigma(u)$, both genetic ornamentations by Lemma~\ref{lem:recording-ideals}.
By the definition of~$\preceq_{\calO}$, they are the meet and the join of~$\rho$ and~$\sigma$ in~$\calO(\sfT)$.
By Theorem~\ref{thm:genetic}, the spine is distributive, and its rank function is~$\rho\mapsto|\Phi(\rho)|=h(\rho)$.
\end{proof}

An element of a finite lattice is \emph{join-irreducible} if it is not the minimum element and has a unique lower cover.
For~$(u,v)\in\calX(\sfT)$, write~$\dn{(u,v)}\coloneqq\{(u',v')\in\calX(\sfT)\colon u'\preceq_\sfT u,\ v'\preceq_\sfT v\}$ for its principal order ideal.

\begin{corollary}\label{cor:join-irreducibles}
The map~$(u,v)\mapsto\Psi(\dn{(u,v)})$ is an order isomorphism from~$\calX(\sfT)$ onto the poset of join-irreducible elements of~$\Spine(\calO(\sfT))$.
The element~$\Psi(\dn{(u,v)})$ is the smallest genetic ornamentation whose ornament at~$u$ contains~$v$, namely
\begin{equation}\label{eq:join-irreducible}
    \Psi(\dn{(u,v)})(z)=
    \begin{cases}
        [z,v]_\sfT & \text{if }z\preceq_\sfT u\,,\\
        \{z\} & \text{otherwise}\,.
    \end{cases}
\end{equation}
\end{corollary}

\begin{proof}
The join-irreducibles of~$J(\calX(\sfT))$ are precisely the principal ideals of~$\calX(\sfT)$~\cite[Proposition~3.4.2]{StanleyEC1}, so the isomorphism follows from Theorem~\ref{thm:genetic}.
Since~$\dn{(u,v)}$ is the smallest ideal containing~$(u,v)$,~$\Psi(\dn{(u,v)})$ is the smallest genetic ornamentation whose ornament at~$u$ contains~$v$.
If~$z\preceq_\sfT u$, the pairs in~$\dn{(u,v)}$ with first coordinate~$z$ are exactly those~$(z,y)$ with~$z\prec_\sfT y\preceq_\sfT v$.
Otherwise there are none.
Applying~$\Psi$ gives~\eqref{eq:join-irreducible}.
\end{proof}

\subsection{Transfer systems}\label{sec:transfer}
We realize~$\calO(\sfT)$ as an interval in a lattice of transfer systems.
Luo and Rognerud's criterion for the spine~\cite[Proposition~7.3]{LuoRognerud26} then gives a second proof of Lemma~\ref{lem:maxchain-implies-genetic}.

\begin{definition}[{\cite[Definition~2.5]{LuoRognerud26}}]\label{def:transfer}
Let~$(L,\preceq)$ be a finite lattice.
A \emph{transfer system} on~$L$ is a subset~$R\subseteq\{(x,y)\in L\times L\colon x\preceq y\}$ such that, for all~$x,y,z\in L$:
\begin{enumerate}
    \renewcommand{\theenumi}{\roman{enumi}}
    \item\label{ts:refl} $(x,x)\in R$;
    \item\label{ts:trans} if~$(x,y)\in R$ and~$(y,z)\in R$, then~$(x,z)\in R$;
    \item\label{ts:restr} if~$(x,y)\in R$ and~$z\preceq y$, then~$(x\wedge z,z)\in R$.
\end{enumerate}
\end{definition}

By~(\ref{ts:refl}),~(\ref{ts:trans}), and the antisymmetry of~$\preceq$, every transfer system on~$L$ is a partial order on~$L$.
The poset~$\Trs(L)$ of transfer systems on~$L$, ordered by inclusion, is a lattice~\cite[Theorem~2.6 and Proposition~3.4]{LuoRognerud26}.

In~$(V,\preceq_\sfT)$, incomparable vertices have no common upper bound, so this poset is a lattice only if~$\sfT$ is a path.
Let~$(\sfT^{+},\preceq_{\sfT^{+}})$ be the poset obtained from~$(V,\preceq_\sfT)$ by adjoining a maximum~$\infty$.
It is a lattice in which the meet of two vertices of~$\sfT$ is their deepest common ancestor.
Let~$\Delta\coloneqq\{(v,v)\colon v\in\sfT^{+}\}$.
For~$\rho\in\calO(\sfT)$, let~$\Theta(\rho)\coloneqq\Delta\cup\Phi(\rho)$; by anchoring,~$(u,v)\in\Theta(\rho)$ if and only if~$v\in\rho(u)$, for all~$u,v\in V$.

\begin{proposition}\label{prop:transfer}
The map~$\Theta$ is an order isomorphism
\[
    \Theta\colon\calO(\sfT)\longrightarrow[\Delta,\Delta\cup\calX(\sfT)]_{\Trs(\sfT^{+})}\,.
\]
\end{proposition}

\begin{proof}
Let~$\rho\in\calO(\sfT)$.
Condition~(\ref{ts:refl}) holds for~$\Theta(\rho)$ because~$\Delta\subseteq\Theta(\rho)$, and~(\ref{ts:trans}) and~(\ref{ts:restr}) are immediate for pairs in~$\Delta$.
So for~(\ref{ts:trans}), consider~$(u,w),(w,y)\in\Phi(\rho)$.
Then~$\rho(u)$ and~$\rho(w)$ both contain~$w$, so~$\rho(w)\subseteq\rho(u)$ by~\eqref{eq:nested-or-disjoint}, and hence~$(u,y)\in\Phi(\rho)$.
For~(\ref{ts:restr}), consider~$(u,v)\in\Phi(\rho)$ and~$z\preceq_{\sfT^{+}}v$.
Since~$u$ and~$z$ are both ancestors of~$v$, they are comparable.
If~$z\preceq_\sfT u$, then~$(u\wedge z,z)=(z,z)\in\Delta$; if~$u\prec_\sfT z$, then~$(u\wedge z,z)=(u,z)\in\Phi(\rho)$, because~$\rho(u)$ is connected and~$z\in[u,v]_\sfT$.
Thus~$\Theta(\rho)$ is a transfer system with~$\Delta=\Theta(\hat0)\subseteq\Theta(\rho)\subseteq\Theta(\hat1)=\Delta\cup\calX(\sfT)$.

Conversely, let~$R$ be a transfer system on~$\sfT^{+}$ contained in~$\Delta\cup\calX(\sfT)$, and let~$\rho\coloneqq\Psi(R\setminus\Delta)$, so that~$\rho(u)=\{v\in V\colon(u,v)\in R\}$ for each~$u\in V$.
If~$v\in\rho(u)$ and~$z\in[u,v]_\sfT$, then~(\ref{ts:restr}) gives~$(u\wedge z,z)=(u,z)\in R$; so~$\rho(u)$ is connected.
If~$x\in\rho(u)\cap\rho(w)$, then~$u$ and~$w$ are comparable because both are ancestors of~$x$; by symmetry, suppose~$u\preceq_\sfT w$.
Then~$w\in[u,x]_\sfT\subseteq\rho(u)$, and~(\ref{ts:trans}) gives~$(u,y)\in R$ for every~$y\in\rho(w)$; hence~$\rho(w)\subseteq\rho(u)$.
Thus~$\rho$ is an ornamentation, and~$\Theta(\rho)=R$ by~\eqref{eq:recording-reconstruction}.
Finally, by~\eqref{eq:recording-reconstruction}, for all~$\rho,\sigma\in\calO(\sfT)$ we have~$\rho\preceq_{\calO}\sigma$ if and only if~$\Theta(\rho)\subseteq\Theta(\sigma)$.
\end{proof}

For~$\sfT=\sfP_n$, deleting~$(\infty,\infty)$ identifies~$[\Delta,\Delta\cup\calX(\sfT)]_{\Trs(\sfT^{+})}$ with~$\Trs(\sfP_n)$, since no other pair in~$\Delta\cup\calX(\sfT)$ involves~$\infty$.
Since~$\calO(\sfP_n)$ is isomorphic to a Tamari lattice, this recovers the result of Balchin, Barnes, and Roitzheim~\cite[Theorem~25]{BBR21} that the lattice of transfer systems on a finite chain is isomorphic to a Tamari lattice (see also~\cite[Theorem~10.2]{LuoRognerud26}).

\begin{proposition}[{\cite[Theorem~2.6 and Proposition~7.3]{LuoRognerud26}}]\label{prop:transfer-spine}
Let~$R$ be a transfer system on a finite lattice~$(L,\preceq)$.
Then~$R\in\Spine(\Trs(L))$ if and only if~$(x',y')\in R$ for all~$(x,y)\in R$ and all~$x',y'\in L$ such that~$x'\preceq y'$,~$x'\preceq x$,~$y'\preceq y$, and~$y'\not\preceq x$.
\end{proposition}

\begin{proof}[Alternative proof of Lemma~\ref{lem:maxchain-implies-genetic}]
The transfer system~$\Delta\cup\calX(\sfT)$ satisfies the criterion of Proposition~\ref{prop:transfer-spine}: if~$(x,y)\in \Delta\cup\calX(\sfT)$ and~$x',y'\in\sfT^{+}$ satisfy~$x'\preceq_{\sfT^{+}}y'$,~$x'\preceq_{\sfT^{+}}x$,~$y'\preceq_{\sfT^{+}}y$, and~$y'\not\preceq_{\sfT^{+}}x$, then~$x\ne y$, so~$y\in V$; hence~$x',y'\in V$ and~$(x',y')\in \Delta\cup\calX(\sfT)$.
So~$\Delta\cup\calX(\sfT)$ lies on a chain of maximum length in~$\Trs(\sfT^{+})$.
Let~$\rho\in\Spine(\calO(\sfT))$.
By Proposition~\ref{prop:transfer},~$\Theta(\rho)$ lies on a chain of maximum length in~$[\Delta,\Delta\cup\calX(\sfT)]_{\Trs(\sfT^{+})}$, so~$\Theta(\rho)\in\Spine(\Trs(\sfT^{+}))$.
Now let~$w$ be a nonroot vertex with parent~$u$, and suppose that~$\rho(w)\ne\{w\}$.
Take~$v\in\rho(w)\setminus\{w\}$.
Since~$(w,v)\in\Theta(\rho)$ and~$u\prec_\sfT w\prec_\sfT v$, Proposition~\ref{prop:transfer-spine} gives~$(u,v)\in\Theta(\rho)$, that is,~$v\in\rho(u)$.
Thus~$\rho(u)$ and~$\rho(w)$ both contain~$v$, so~$\rho(w)\subseteq\rho(u)$ by~\eqref{eq:nested-or-disjoint}.
Hence~$\rho$ is genetic.
\end{proof}

By Proposition~\ref{prop:transfer} and Theorem~\ref{thm:genetic},~$\Spine([\Delta,\Delta\cup\calX(\sfT)]_{\Trs(\sfT^{+})})$ consists of the transfer systems~$\Delta\cup I$ with~$I\in J(\calX(\sfT))$.

\section{Decomposing the poset of pairs}\label{sec:enumeration}

For a finite poset~$P$, write~$\LE(P)$ for its number of linear extensions.
By Corollary~\ref{cor:linext}, the number~$\ML(\sfT)$ of chains of maximum length in~$\calO(\sfT)$ equals~$\LE(\calX(\sfT))$.
To compute~$\LE(\calX(\sfT))$, we first decompose~$\calX(\sfT)$ at the root.

If~$P$ is a disjoint union~$P_1\sqcup\cdots\sqcup P_k$ of finite posets, then, by interleaving linear extensions of~$P_1,\dots,P_k$, we obtain~\cite[Example~3.5.4]{StanleyEC1}
\begin{equation}\label{eq:LE-disjoint}
    \LE(P)=\binom{|P_1|+\cdots+|P_k|}{|P_1|,\dots,|P_k|}\prod_{i=1}^{k}\LE(P_i)\,.
\end{equation}
We use the conventions that an empty sum equals~$0$, that an empty product and an empty multinomial coefficient equal~$1$, and that the empty poset has one linear extension.

\begin{definition}\label{def:branch}
    Let~$c_1,\dots,c_k$ be the children of the root~$r$.
    For each~$i\in[k]$, let~$\sfT_i$ be the subtree induced by~$\{r\}\cup\desc{c_i}$, rooted at~$r$.
\end{definition}

\begin{proposition}\label{prop:branch}
    Let~$\sfT_1,\dots,\sfT_k$ be the trees of Definition~\ref{def:branch}.
    Then~$\calX(\sfT)=\bigsqcup_{i=1}^{k}\calX(\sfT_i)$ as posets.
    Hence~$H(\sfT)=\sum_{i=1}^{k}H(\sfT_i)$ and
    \begin{equation}\label{eq:branch-ML}
        \ML(\sfT)=H(\sfT)!\prod_{i=1}^{k}\frac{\ML(\sfT_i)}{H(\sfT_i)!}\,.
    \end{equation}
\end{proposition}

\begin{proof}
    For~$(u,v)\in\calX(\sfT)$, the vertex~$v$ belongs to exactly one tree~$\sfT_i$, which also contains all of~$v$'s ancestors.
    Thus~$(u,v)$ belongs to exactly one~$\calX(\sfT_i)$.
    Pairs in~$\calX(\sfT_i)$ and~$\calX(\sfT_j)$ with~$i\ne j$ are incomparable because their second coordinates are incomparable.
    The formulas follow from~$H(\sfT)=|\calX(\sfT)|$, Corollary~\ref{cor:linext}, and~\eqref{eq:LE-disjoint}.
\end{proof}

By Proposition~\ref{prop:branch}, it suffices to count chains of maximum length when the root has a unique child.

\begin{lemma}\label{lem:columns}
    The map~$(u,v)\mapsto(\dep(u),v)$ is an order isomorphism from~${(\calX(\sfT),\preceq_{\calX})}$ onto the set~$\{(i,v)\in\mathbb{Z}_{\ge0}\times V\colon i<\dep(v)\}$, ordered by~${(i,v)\preceq(i',v')}$ if~$i\le i'$ and~$v\preceq_\sfT v'$.
\end{lemma}

\begin{proof}
    A vertex~$v\neq r$ has exactly one proper ancestor of depth~$i$ for each~$0\le i\le\dep(v)-1$, so the map is a bijection.
    Both orders require~$v\preceq_\sfT v'$, so it remains to compare first coordinates.
    Let~$u\prec_\sfT v\preceq_\sfT v'$ and~$u'\prec_\sfT v'$.
    The vertices~$u$,~$u'$, and~$v$ lie on the path from~$r$ to~$v'$, where~$\preceq_\sfT$ agrees with the order of depths.
    Thus~$u\preceq_\sfT u'$ if and only if~${\dep(u)\le\dep(u')}$, so the map preserves and reflects order.
\end{proof}

The \emph{column coordinates} of a pair~$(u,v)$ are~$(\dep(u),v)$; in particular,~$(0,v)$ represents~$(r,v)$.
We retain the order symbol~$\preceq_{\calX}$ for pairs in column coordinates.

\begin{definition}\label{def:column}
    For a nonroot vertex~$v$, its \emph{column} is the chain
    \[
        C_v\coloneqq \{(i,v)\colon 0\le i\le\dep(v)-1\}\,.
    \]
\end{definition}

\begin{figure}[ht]
\centering
\includegraphics{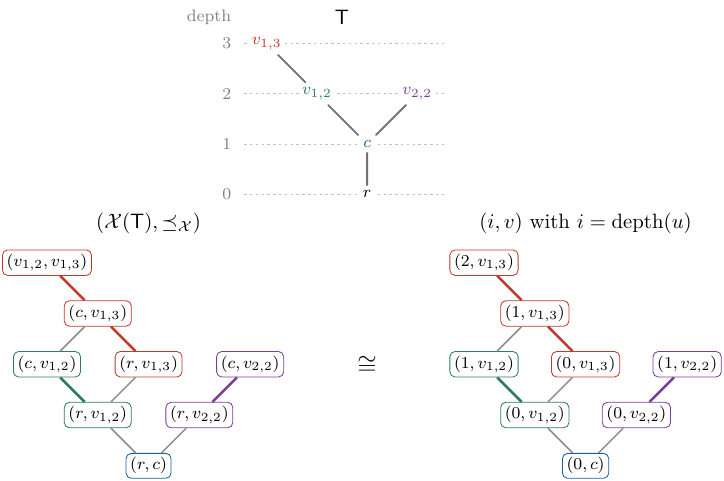}
\caption{The elements of~$\calX(\sfT)$ written as pairs of vertices~$(u,v)$ (left) and in column coordinates~$(\dep(u),v)$ (right).
Colors indicate the second coordinate.}
\label{fig:columns}
\end{figure}

In Figure~\ref{fig:columns}, the colored edges lie within a column.
The gray edges join the columns of a vertex and its parent.
In our figures, a vertex named~$v_{j,t}$ has depth~$t$, and vertices with the same index~$j$ lie on a common path from~$r$.

\begin{definition}\label{def:block}
    The \emph{block} of a nonroot vertex~$v$ is the induced subposet of pairs with second coordinate in~$\desc{v}$, namely
    \[
        \calB_v\coloneqq \bigcup_{x\in\desc{v}}C_x
        \subseteq\calX(\sfT)\,.
    \]
\end{definition}

\begin{samepage}
\begin{lemma}\label{lem:branch-skeleton}
    Let~$\sfT$ be a rooted tree whose root has exactly one child~$c$.
    Then~$(0,c)$ is the minimum of~$\calX(\sfT)$, and
    \[
        \calX(\sfT)\setminus\{(0,c)\}=\bigsqcup_{\text{$w$ child of $c$}}\calB_w
        \qquad\text{as posets}\,.
    \]
    Hence~$H(\sfT)=1+\sum_{\text{$w$ child of $c$}}|\calB_w|$ and
    \begin{equation}\label{eq:skeleton}
        \ML(\sfT)=\bigl(H(\sfT)-1\bigr)!\prod_{\text{$w$ child of $c$}}\frac{\LE(\calB_w)}{|\calB_w|!}\,.
    \end{equation}
\end{lemma}
\end{samepage}

\begin{proof}
    By Lemma~\ref{lem:columns},~$C_c=\{(0,c)\}$ and~$(0,c)$ lies below every other element, so every linear extension of~$\calX(\sfT)$ starts with~$(0,c)$.
    Every column other than~$C_c$ lies in exactly one block~$\calB_w$, and columns in distinct blocks have incomparable second coordinates.
    Then~\eqref{eq:skeleton} follows from~\eqref{eq:LE-disjoint}, applied to the disjoint union of the blocks~$\calB_w$.
\end{proof}

\begin{definition}[{\cite[Definitions~2.3 and~2.14]{AdinRoichman15}}]\label{def:shifted-diagram}
    Let~$\lambda=(\lambda_1>\lambda_2>\cdots>\lambda_\ell>0)$ be a strict partition with~$\ell$ parts.
    Its \emph{size}~$|\lambda|$ is the sum of its parts.
    Its \emph{shifted diagram} is the poset
    \[
        \calD_\lambda\coloneqq \{(a,b)\in\mathbb{Z}_{>0}^2\colon 1\le a\le \ell,\ a\le b\le a+\lambda_a-1\}\,,
    \]
    ordered componentwise.
    For~$m\ge0$, set~$\delta_m\coloneqq (m,m-1,\dots,1)$, with~$\delta_0\coloneqq \varnothing$.
    We call~$\calD_{\delta_m}$ a \emph{shifted staircase}.
    A \emph{standard shifted Young tableau} of shape~$\lambda$ is a labeling of~$\calD_\lambda$ by~$1,\ldots,|\lambda|$ that is increasing for the componentwise order.
    The corresponding linear extension is the list of cells in increasing label order.
    We write~$g^\lambda\coloneqq\LE(\calD_\lambda)$ for the number of these tableaux.
\end{definition}

We call the rows and columns of a shifted diagram its \emph{diagram rows} and \emph{diagram columns}, to distinguish them from the columns~$C_v$ of Definition~\ref{def:column}.

\begin{figure}[ht]
\centering
\includegraphics{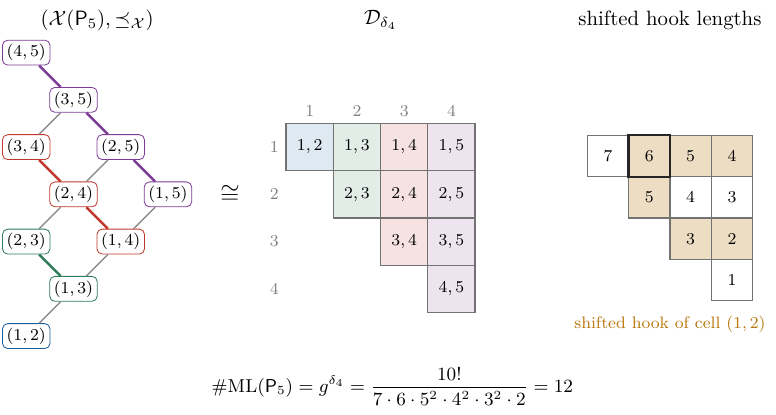}
\caption{The pair~$(u,v)$ in~$\calX(\sfP_5)$ (left) corresponds to the cell~$(u,v-1)$ of~$\calD_{\delta_4}$ (middle).
On the right are the shifted hook lengths of~$\calD_{\delta_4}$, with the shifted hook of cell~$(1,2)$, which holds the pair~$(1,3)$, shaded.}
\label{fig:staircase}
\end{figure}

The \emph{shifted hook}~\cite[Definition~5.8]{AdinRoichman15} of a cell~$(a,b)$ in~$\calD_{\delta_m}$ consists of the cell itself, the cells to its right in diagram row~$a$, and the cells below it in diagram column~$b$.
If~$b<m$, the hook also includes all of diagram row~$b+1$.
The shifted hook has~$2m+1-a-b$ cells (Figure~\ref{fig:staircase}).

\begin{proposition}[{cf.~\cite[Corollary~3.4]{FishelNelson2014}}]\label{prop:paths}
    For~$n\ge1$, the map sending a pair of vertices~$(u,v)$ to~$(u,v-1)$ (in column coordinates,~$(i,v)\mapsto(i+1,\dep(v))$) is an order isomorphism from~$\calX(\sfP_n)$ onto the shifted staircase~$\calD_{\delta_{n-1}}$; hence~$H(\sfP_n)=\binom{n}{2}$ and
    \[
        \ML(\sfP_n)=g^{\delta_{n-1}}=\binom{n}{2}!\Bigm/\prod_{1\le a\le b\le n-1}(a+b-1)\,.
    \]
\end{proposition}

\begin{proof}
    The map is a bijection that preserves and reflects the componentwise order (Figure~\ref{fig:staircase}).
    Counting cells gives~$H(\sfP_n)$ by Corollary~\ref{cor:chain-length}.
    For~$\ML(\sfP_n)$, apply Corollary~\ref{cor:linext} and the shifted hook-length formula~\cite[Theorem~5.9]{AdinRoichman15}, reindexing the hook lengths~$2n-1-a-b$ by~$(a,b)\mapsto(n-b,n-a)$.
\end{proof}

By Proposition~\ref{prop:paths} and Corollary~\ref{cor:linext}, we recover the ordinary Tamari-lattice case of the theorem of Fishel and Nelson~\cite[Theorem~3.3 and Corollary~3.4]{FishelNelson2014}.
The maximum length~$\binom{n}{2}$ was previously determined by Markowsky~\cite{Markowsky92}.
Fishel and Nelson also identified the elements on chains of maximum length in the Tamari lattice with the strict partitions whose shifted diagrams are contained in~$\calD_{\delta_{n-1}}$~\cite[Section~3]{FishelNelson2014}.

\begin{samepage}
\begin{definition}\label{def:spider}
    For~$k\ge0$ and~$n_1,\dots,n_k\ge1$, the \emph{spider}~$\sfSp(n_1,\dots,n_k)$ is obtained by identifying the roots of the paths~$\sfP_{n_1+1},\dots,\sfP_{n_k+1}$.
    The common root is the \emph{center}.
    The~$j$th \emph{leg} is the set of~$n_j$ noncentral vertices on the~$j$th path.
    The spider with no legs is the one-vertex tree.
\end{definition}
\end{samepage}

The set~$\desc{w}$ induces a spider centered at~$w$ if and only if every proper descendant of~$w$ has at most one child; in that case, we call the legs of this spider the \emph{legs of~$w$}.
Let~$w$ be such a vertex and set~$d\coloneqq\dep(w)$.
On a leg with~$n$ vertices beginning at a child~$w'$ of~$w$, the vertices have depths~$d+1,\ldots,d+n$.
The block~$\calB_{w'}$ is the union of their columns, so
\begin{equation}\label{eq:leg-size}
    |\calB_{w'}|=\sum_{t=1}^{n}(t+d)=\frac{n(n+2d+1)}{2}\,.
\end{equation}

\begin{samepage}
\begin{lemma}\label{lem:leg-anatomy}
    Let~$w$ be a vertex of depth~$d\ge1$ such that~$\desc{w}$ induces a spider centered at~$w$.
    Fix a leg of the spider, and let~$n$ be the number of vertices of the leg.
    \begin{enumerate}
        \item The subposet induced by~$C_w$ and the columns of the chosen leg is isomorphic to the subposet~${\calD_{\delta_{n+d}}\setminus\calD_{\delta_{d-1}}}$ of~$\calD_{\delta_{n+d}}$.
        The isomorphism sends~$(i,v)$ to~${(i+1,\dep(v))}$.
        The column~$C_w$ maps onto diagram column~$d$, and the columns of the chosen leg onto diagram columns~$d+1,\dots,n+d$.
        \item Elements in the columns of distinct legs are incomparable.
    \end{enumerate}
\end{lemma}

\begin{proof}
    Since~$w$ and the vertices of the chosen leg lie on a common path starting at~$r$, on which~$\preceq_\sfT$ agrees with the order of depths, Lemma~\ref{lem:columns} shows that~$(i,v)\mapsto(i+1,\dep(v))$ is an order isomorphism onto its image.
    The column of a vertex of depth~$b$ maps onto the full diagram column~$b$, so the image consists of diagram columns~$d,d+1,\dots,n+d$.
    The remaining diagram columns~$1,\dots,d-1$ of~$\calD_{\delta_{n+d}}$ form~$\calD_{\delta_{d-1}}$, proving~(1).
    Vertices on distinct legs are incomparable, so~(2) also follows from Lemma~\ref{lem:columns}.
\end{proof}
\end{samepage}

\section{Trees of branching depth at most 1}\label{sec:bd-one}

\begin{definition}\label{def:branching-depth}
    The \emph{branching depth} of~$\sfT$ is
    \[
        \bd(\sfT)\coloneqq \max\bigl\{\dep(v)\colon\text{$v\in V$ has at least two children}\bigr\}\,.
    \]
    We set~$\bd(\sfT)\coloneqq 0$ when every vertex of~$\sfT$ has at most one child.
\end{definition}

\begin{definition}\label{def:stemmed}
    For~$k\ge0$ and~$n_1,\dots,n_k\ge1$, the \emph{stemmed spider}~$\sfSp^{+}(n_1,\dots,n_k)$ is obtained from~$\sfSp(n_1,\dots,n_k)$ by adjoining a new root below the center.
    Its center and legs are those of the spider.
\end{definition}

In particular, the stemmed spider with no legs is~$\sfP_2$, and~$\sfSp^{+}(n)=\sfP_{n+2}$.
The tree of Figure~\ref{fig:columns} is the stemmed spider~$\sfSp^{+}(2,1)$; its center is the vertex~$c$ of depth~$1$.

\begin{observation}\label{obs:bd-one-branches}
    We have~$\bd(\sfT)\le1$ if and only if every vertex of depth at least~$2$ has at most one child, that is, if and only if each tree~$\sfT_i$ of Definition~\ref{def:branch} is a stemmed spider.
\end{observation}

Spiders (in particular paths and stars) and trees of depth at most~$2$ have branching depth at most~$1$.

\begin{proposition}\label{prop:stemmed-spider}
    Let~$\sfT=\sfSp^{+}(n_1,\dots,n_k)$.
    Then~$H(\sfT)=1+\sum_{j=1}^{k}\bigl(\tbinom{n_j+2}{2}-1\bigr)$ and
    \begin{equation}\label{eq:stemmed}
        \ML(\sfT)=\bigl(H(\sfT)-1\bigr)!\;\prod_{j=1}^{k}\frac{\binom{n_j+2}{2}}{\prod_{1\le a\le b\le n_j+1}(a+b-1)}\,.
    \end{equation}
\end{proposition}

\begin{proof}
    By Lemma~\ref{lem:leg-anatomy}(1) with~$d=1$, the block of the first vertex of the~$j$th leg is isomorphic to~${\calD_{\delta_{n_j+1}}\setminus\{(1,1)\}}$.
    The block has~$\binom{n_j+2}{2}-1$ elements and~$g^{\delta_{n_j+1}}$ linear extensions, since~$(1,1)$ comes first in every linear extension of~$\calD_{\delta_{n_j+1}}$.
    Substitute the value of~$g^{\delta_{n_j+1}}$ from Proposition~\ref{prop:paths} into~\eqref{eq:skeleton}.
\end{proof}

\begin{figure}[ht]
\centering
\includegraphics{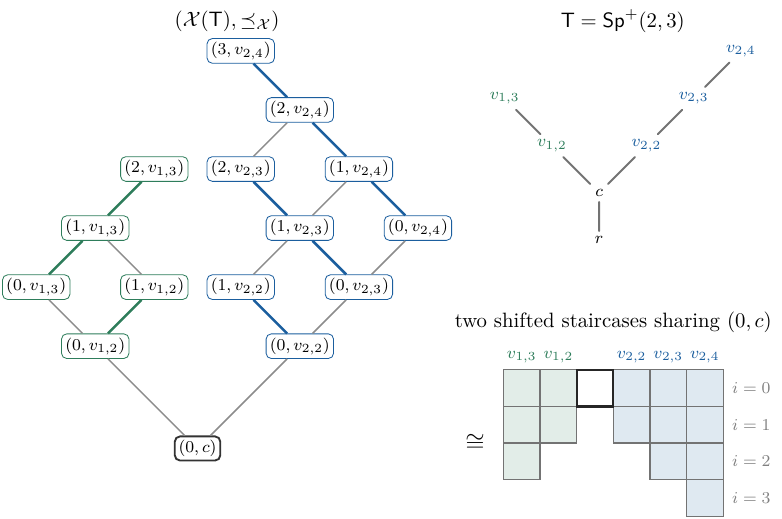}
\caption{The poset~$\calX(\sfSp^{+}(2,3))$ as two shifted staircases sharing their minimum~$(0,c)$, with one staircase reflected.
Colors indicate the leg containing the second coordinate.}
\label{fig:glued-staircases}
\end{figure}

Figure~\ref{fig:glued-staircases} shows the case~$\sfSp^{+}(2,3)$.
In the staircase drawing, diagram columns correspond to the columns~$C_c$ and~$C_{v_{j,t}}$, and diagram rows have constant~$i$.

\begin{corollary}\label{cor:bd-one}
    Let~$\sfT$ be a rooted tree with~$\bd(\sfT)\le1$, and let~$\sfT_1,\dots,\sfT_k$ be as in Definition~\ref{def:branch}.
    For each~$i\in[k]$, let~$n_{i,1},\dots,n_{i,k_i}$ be the numbers of vertices in the legs of the stemmed spider~$\sfT_i$.
    Then
    \[
        \ML(\sfT)=H(\sfT)!\,\prod_{i=1}^{k}\Biggl(
        \frac{1}{H(\sfT_i)}
        \prod_{j=1}^{k_i}\frac{\binom{n_{i,j}+2}{2}}{\prod_{1\le a\le b\le n_{i,j}+1}(a+b-1)}\Biggr)\,.
    \]
\end{corollary}

\begin{proof}
    By Observation~\ref{obs:bd-one-branches}, each~$\sfT_i$ is a stemmed spider; substitute~\eqref{eq:stemmed} for~$\ML(\sfT_i)$ in~\eqref{eq:branch-ML}.
\end{proof}

\begin{definition}\label{def:depth-two-family}
    For~$k\ge0$ and~$a_1,\dots,a_k\ge0$, let~$\sfD(a_1,\dots,a_k)$ denote the rooted tree with root~$r$, children~$c_1,\dots,c_k$, and~$a_i$ leaves attached to~$c_i$.
\end{definition}

The trees~$\sfD(a_1,\dots,a_k)$ are exactly the rooted trees of depth at most~$2$.

\begin{corollary}\label{cor:depth-two}
    Let~$\sfT=\sfD(a_1,\dots,a_k)$.
    Then
    \begin{equation}\label{eq:depth-two}
        \ML(\sfT)=\bigl(k+2{\textstyle\sum_i a_i}\bigr)!\Bigm/\prod_{i=1}^{k}\bigl((2a_i+1)\,2^{a_i}\bigr)\,.
    \end{equation}
\end{corollary}

\begin{proof}
    The tree~$\sfT_i$ of Definition~\ref{def:branch} is the stemmed spider~$\sfSp^{+}(1,\dots,1)$ with~$a_i$ legs, so Proposition~\ref{prop:stemmed-spider} gives~$H(\sfT_i)=2a_i+1$ and~$\ML(\sfT_i)=(2a_i)!/2^{a_i}$.
    Substituting these values into~\eqref{eq:branch-ML}, with~$H(\sfT)=k+2\sum_ia_i$, yields~\eqref{eq:depth-two}.
\end{proof}

For~$\sfY=\sfD(2)$, this gives~$5!/(5\cdot2^{2})=6$ chains of maximum length, as in Figure~\ref{fig:tree-and-pairposet}.

When~$\sfT$ has depth at most~$2$, each element of~$\calX(\sfT)$ covers at most one element.
The Hasse diagram of~$\calX(\sfT)$ is therefore a rooted forest.
Corollary~\ref{cor:depth-two} also follows from Knuth's tree hook-length formula~\cite[Section~5.1.4, Exercise~20]{Knuth98}, applied to each component, and~\eqref{eq:LE-disjoint}.

\section{Trees of branching depth at most 2}\label{sec:bd-two}

If~$\bd(\sfT)\le2$, then by Proposition~\ref{prop:branch} and Lemma~\ref{lem:branch-skeleton} it suffices to count the linear extensions of~$\calB_w$ for each vertex~$w$ of depth~$2$.
Until Theorem~\ref{thm:bd-two}, fix a vertex~$w$ of depth~$2$ such that~$\desc{w}$ induces a spider centered at~$w$.

\begin{lemma}\label{lem:gap}
    Every linear extension of~$\calB_w$ starts with~$(0,w)$.
    If~$v_1\prec_\sfT\cdots\prec_\sfT v_n$ are the vertices of a leg of~$w$, then in each linear extension of~$\calB_w$, the pairs of~$\calB_{v_1}$ placed between~$(0,w)$ and~$(1,w)$ are~$(0,v_1),\dots,(0,v_q)$ for some~$0\le q\le n$.
\end{lemma}

\begin{figure}[ht]
\centering
\includegraphics{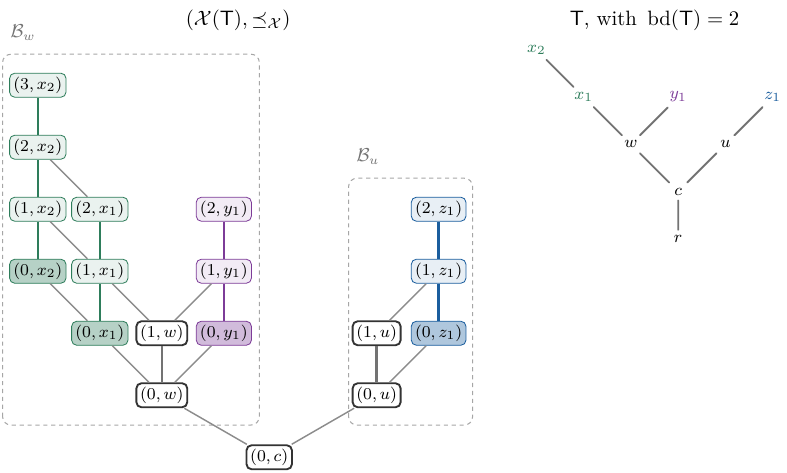}
\caption{Deleting~$(0,c)$ separates~$\calX(\sfT)$ into the blocks~$\calB_w$ and~$\calB_u$.
The darker boxes mark the pairs that linear extensions of~$\calB_w$ and~$\calB_u$ can place between~$(0,w)$ and~$(1,w)$, and between~$(0,u)$ and~$(1,u)$, respectively.
The thick black borders mark the columns~$C_c$,~$C_w$, and~$C_u$.
Colors indicate the leg containing the second coordinate.}
\label{fig:one-gap}
\end{figure}

\begin{proof}
    By Lemma~\ref{lem:columns},~$(0,w)$ is the minimum of~$\calB_w$.
    For~$i\ge1$ and a vertex~$v$ on a leg of~$w$, we have~$(1,w)\preceq_{\calX}(i,v)$, so~$(i,v)$ comes after~$(1,w)$.
    Since the set of elements placed before~$(1,w)$ is an ideal, its intersection with~$\calB_{v_1}$ is~$\{(0,v_1),\dots,(0,v_q)\}$ for some~$q$ (Figure~\ref{fig:one-gap}).
\end{proof}

Fix a leg of~$w$ with vertices~$v_1\prec_\sfT\cdots\prec_\sfT v_n$.
By Lemma~\ref{lem:leg-anatomy}(1),~${(i,v)\mapsto(i+1,\dep(v))}$ is an order isomorphism from~$C_w\cup\calB_{v_1}$ onto~$\calD_{\delta_{n+2}}\setminus\{(1,1)\}$.
If a linear extension of~$\calB_w$ places exactly the pairs~$(0,v_1),\dots,(0,v_q)$ of~$\calB_{v_1}$ between~$(0,w)$ and~$(1,w)$ (Lemma~\ref{lem:gap}), then~$(0,w)$,~$(0,v_1),\dots,(0,v_q)$, and~$(1,w)$ map onto~$\calD_{(q+2,1)}\setminus\{(1,1)\}$, so the pairs of~$\calB_{v_1}$ placed after~$(1,w)$ map onto~$\calD_{\delta_{n+2}}\setminus\calD_{(q+2,1)}$.

\begin{samepage}
\begin{lemma}[{\cite[Lemma~4.2]{AdinKingRoichman11}}]\label{lem:staircase-complement}
    Let~$m\ge0$ and let~$\mu$ be a strict partition with~$\calD_\mu\subseteq\calD_{\delta_m}$.
    Then the complement~$\calD_{\delta_m}\setminus\calD_\mu$, as a subposet of~$\calD_{\delta_m}$, satisfies
    \[
        \LE\bigl(\calD_{\delta_m}\setminus\calD_\mu\bigr)=g^{[m]\setminus\mu}\,,
    \]
    where~$[m]\setminus\mu$ is the strict partition formed by the integers in~$[m]$ absent from~$\mu$.
\end{lemma}
\end{samepage}

\begin{figure}[ht]
\centering
\includegraphics{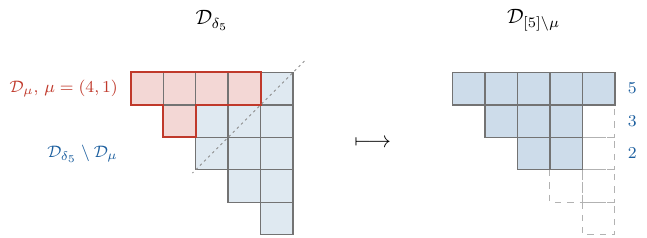}
\caption{The antidiagonal reflection of~$\calD_{\delta_5}\setminus\calD_{(4,1)}$ onto~$\calD_{(5,3,2)}$, with the reflection axis dotted.}
\label{fig:staircase-complement}
\end{figure}

Lemma~\ref{lem:staircase-complement} holds because the antidiagonal reflection~$(a,b)\mapsto(m+1-b,m+1-a)$ reverses the order of~$\calD_{\delta_m}$ and carries~$\calD_{\delta_m}\setminus\calD_\mu$ onto~$\calD_{[m]\setminus\mu}$.
Figure~\ref{fig:staircase-complement} shows the case~$m=5$,~$\mu=(4,1)$.

By Lemma~\ref{lem:staircase-complement},
\begin{equation}\label{eq:leg-completions}
    \LE\bigl(\calD_{\delta_{n+2}}\setminus\calD_{(q+2,1)}\bigr)=g^{[n+2]\setminus(q+2,1)}\,.
\end{equation}
For~$n=3$ and~$q=2$, we have~$[5]\setminus(4,1)=(5,3,2)$ (Figure~\ref{fig:staircase-complement}) and~${g^{(5,3,2)}=54}$.

By Schur's product formula~\cite{Schur1911} (see also Thrall~\cite{Thrall52} and Adin and Roichman~\cite[Theorem~5.7]{AdinRoichman15}), every strict partition~$\lambda$ satisfies
\begin{equation}\label{eq:schur}
    g^{\lambda}
    =\frac{|\lambda|!}{\prod_{a}\lambda_a!}
    \prod_{a<a'}\frac{\lambda_a-\lambda_{a'}}{\lambda_a+\lambda_{a'}}\,.
\end{equation}

\begin{lemma}\label{lem:Anq}
    For~$n\ge1$ and~$0\le q\le n$, we have
    \begin{equation}\label{eq:Anq}
        g^{[n+2]\setminus(q+2,1)}
        =\frac{\bigl(\tfrac{n(n+5)}{2}-q\bigr)!}{\prod_{1\le a\le b\le n+2}(a+b-1)}\;
        \frac{(n+2)(n+3)(q+1)}{4(q+3)}\,
        \frac{(n+q+4)!}{(q+2)!\,(n-q)!}\,.
    \end{equation}
\end{lemma}

\begin{proof}
    Let~$\lambda\coloneqq[n+2]\setminus(q+2,1)$, of size~$|\lambda|=\tfrac{n(n+5)}{2}-q$.
    The parts of~$\delta_{n+2}$ are the parts of~$\lambda$ together with~$1$ and~$q+2$.
    Dividing Schur's formula~\eqref{eq:schur} for~$\lambda$ by the formula for~$\delta_{n+2}$ therefore leaves~$|\lambda|!/\tbinom{n+3}{2}!$, the factor~$1!\,(q+2)!$, and one factor~$\tfrac{s+t}{|s-t|}$ for each pair~${\{s,t\}\subseteq[n+2]}$ that contains~$1$ or~$q+2$:
    \[
        \frac{g^\lambda}{g^{\delta_{n+2}}}
        =\frac{|\lambda|!}{\binom{n+3}{2}!}\;1!\,(q+2)!\,
        \prod_{t=2}^{n+2}\frac{t+1}{t-1}\,
        \prod_{\substack{2\le t\le n+2\\ t\ne q+2}}\frac{q+2+t}{|q+2-t|}\,.
    \]
    The two products equal~$\tfrac{(n+2)(n+3)}{2}$ and~$(n+q+4)!\big/\bigl(2(q+2)\,(q+3)!\,q!\,(n-q)!\bigr)$, respectively.
    Substituting these values and~$g^{\delta_{n+2}}=\binom{n+3}{2}!\big/\prod_{1\le a\le b\le n+2}(a+b-1)$ from Proposition~\ref{prop:paths} gives~\eqref{eq:Anq}.
\end{proof}

\begin{definition}\label{def:leg-poly}
    For~$n\ge1$, let
    \begin{equation}\label{eq:leg-poly}
        L_n(z)\coloneqq
        \frac{(n+2)(n+3)}{4\prod_{1\le a\le b\le n+2}(a+b-1)}
        \sum_{q=0}^{n}\frac{(q+1)\,(n+q+4)!}{(q+3)\,q!\,(q+2)!\,(n-q)!}\,z^{q}\,.
    \end{equation}
\end{definition}

By Lemma~\ref{lem:Anq}, for~$0\le q\le n$, the coefficient of~$z^q$ in~$L_n(z)$ is
\begin{equation}\label{eq:leg-coeff}
    \frac{g^{[n+2]\setminus(q+2,1)}}{q!\,\bigl(\tfrac{n(n+5)}{2}-q\bigr)!}\,.
\end{equation}

\begin{samepage}
\begin{proposition}\label{prop:two-chain}
    Let~$w$ be a vertex of depth~$2$ such that~$\desc{w}$ induces a spider centered at~$w$.
    Assume that the~$k$ legs of~$w$ have~$n_1,\dots,n_k$ vertices, respectively.
    Then
    \begin{equation}\label{eq:two-chain}
        \LE(\calB_w)
        =\sum_{q=0}^{\lvert\desc{w}\rvert-1}
        q!\,\bigl(|\calB_w|-2-q\bigr)!\;[z^q]\prod_{j=1}^{k}L_{n_j}(z)\,,
    \end{equation}
    where~$L_n(z)$ is defined by~\eqref{eq:leg-poly} and~$[z^q]$ extracts the coefficient of~$z^q$.
\end{proposition}
\end{samepage}

\begin{proof}
    For each~$j\in[k]$, let~$v^{(j)}_1\prec_\sfT\cdots\prec_\sfT v^{(j)}_{n_j}$ be the vertices of the~$j$th leg.
    Consider a linear extension of~$\calB_w$, and let~$Q$ be the set of pairs that it places between~$(0,w)$ and~$(1,w)$.
    By Lemma~\ref{lem:gap}, the linear extension starts with~$(0,w)$, and~$Q=Q_1\sqcup\cdots\sqcup Q_k$, where~$Q_j$ is the chain~$(0,v^{(j)}_1)\prec_{\calX}\cdots\prec_{\calX}(0,v^{(j)}_{q_j})$ for some~$0\le q_j\le n_j$.
    Fix~$q_1,\dots,q_k$, which determine~$Q$, and let~$R\coloneqq\calB_w\setminus(C_w\cup Q)$.

    Write~$q\coloneqq q_1+\cdots+q_k$,~$R_j\coloneqq R\cap\calB_{v^{(j)}_1}$, and~$\lambda^{(j)}\coloneqq[n_j+2]\setminus(q_j+2,1)$.
    By Lemma~\ref{lem:leg-anatomy}(2),~$Q$ and~$R$ are the disjoint unions of~$Q_1,\dots,Q_k$ and of~$R_1,\dots,R_k$, as posets.
    By Lemma~\ref{lem:leg-anatomy}(1),~\eqref{eq:leg-size}, and~\eqref{eq:leg-completions}, each~$R_j$ has~$\tfrac{n_j(n_j+5)}{2}-q_j$ elements and~$g^{\lambda^{(j)}}$ linear extensions.
    Applying~\eqref{eq:LE-disjoint} to~$Q$ and to~$R$ and summing~$\LE(Q)\,\LE(R)$ over~$q_1+\cdots+q_k=q$ gives the number of linear extensions of~$\calB_w$ that place exactly~$q$ pairs between~$(0,w)$ and~$(1,w)$:
    \[
        q!\,\bigl(|\calB_w|-2-q\bigr)!
        \sum_{\substack{q_1+\cdots+q_k=q\\0\le q_j\le n_j}}
        \prod_{j=1}^{k}\frac{g^{\lambda^{(j)}}}{q_j!\,\bigl(\tfrac{n_j(n_j+5)}{2}-q_j\bigr)!}\,.
    \]
    By~\eqref{eq:leg-coeff}, each factor of the product is the coefficient of~$z^{q_j}$ in~$L_{n_j}(z)$, so the sum is~$[z^q]\prod_jL_{n_j}(z)$.
    Summing over~$q$ proves~\eqref{eq:two-chain}.
\end{proof}

\begin{samepage}
\begin{theorem}\label{thm:bd-two}
    Let~$\sfT$ be a rooted tree with~$\bd(\sfT)\le2$.
    Let~$c_1,\dots,c_k$ and~$\sfT_1,\dots,\sfT_k$ be as in Definition~\ref{def:branch}.
    For each vertex~$w$ of depth~$2$, the set~$\desc{w}$ induces a spider centered at~$w$; let~$k_w\ge0$ be the number of its legs, and let~$n_{w,1},\dots,n_{w,k_w}$ be their numbers of vertices.
    Then
    \begin{equation}\label{eq:bd-two-explicit}
        \ML(\sfT)
        =\frac{H(\sfT)!}{\prod_{i=1}^{k}H(\sfT_i)}
        \prod_{\substack{w\in V\\ \dep(w)=2}}\Biggl(
        \frac{1}{|\calB_w|!}
        \sum_{q=0}^{\lvert\desc{w}\rvert-1}
        q!\,\bigl(|\calB_w|-2-q\bigr)!\;[z^q]\prod_{j=1}^{k_w}L_{n_{w,j}}(z)\Biggr)\,.
    \end{equation}
    Here~$|\calB_w|=2+\sum_{j=1}^{k_w}\tfrac{n_{w,j}(n_{w,j}+5)}{2}$ and~$H(\sfT_i)=1+\sum_{\text{$w$ child of $c_i$}}|\calB_w|$.
\end{theorem}
\end{samepage}

\begin{proof}
    Since~$\bd(\sfT)\le2$, every proper descendant of a vertex~$w$ of depth~$2$ has at most one child, so~$\desc{w}$ induces a spider centered at~$w$.
    Each block~$\calB_w$ is the union of~$C_w$, which has two elements, and the blocks~$\calB_{w'}$ of the children~$w'$ of~$w$, of sizes given by~\eqref{eq:leg-size}.
    Lemma~\ref{lem:branch-skeleton} then gives~$H(\sfT_i)$.
    Applying Proposition~\ref{prop:two-chain} to each block,~\eqref{eq:skeleton} to each~$\sfT_i$, and then~\eqref{eq:branch-ML} gives~\eqref{eq:bd-two-explicit}.
\end{proof}

\begin{definition}\label{def:broom}
    For~$d,k\ge0$, the \emph{broom}~$\sfB_{d,k}$ is obtained from a path~$u_0,u_1,\dots,u_d$, rooted at~$u_0$, by attaching~$k$ leaves to~$u_d$.
\end{definition}

In particular,~$\sfB_{0,k}$ is the star~$\sfS_k$, while~$\sfB_{d,0}=\sfP_{d+1}$ and~$\sfB_{d,1}=\sfP_{d+2}$ are paths.

Defant and Sack~\cite[Sections~1 and~7]{defant2024operahedronlattices} also studied brooms, with a different indexing convention.
Abram, Bastidas, G\'elinas, Pilaud, and Sack~\cite[Section~5.3]{abram2025ornamentation} counted the elements of the ornamentation lattice of a broom.

\begin{corollary}\label{cor:medium-brooms}
    For every~$k\ge0$, we have
    \[
        \ML(\sfB_{2,k})
        =\sum_{q=0}^{k}\frac{k!}{(k-q)!}\,\frac{(3k-q)!}{2^{q}\,6^{k-q}}\,.
    \]
\end{corollary}

\begin{proof}
    By~\eqref{eq:skeleton},~$\ML(\sfB_{2,k})=\LE(\calB_{u_2})$.
    Each of the~$k$ legs of~$u_2$ has one vertex, so substitute~$L_1(z)=\tfrac16+\tfrac12z$ and~$|\calB_{u_2}|=3k+2$ into~\eqref{eq:two-chain} and apply the binomial theorem.
\end{proof}

For~$k=0,\dots,4$, the values of~$\ML(\sfB_{2,k})$ are~$1$,~$2$,~$52$,~$5160$, and~$1{,}283{,}520$.
For a leg with~$n$ vertices, the product~$\prod_{1\le a\le b\le n+1}(a+b-1)$ in Corollary~\ref{cor:bd-one} contains the factors~$n+1$ and~$n+2$, and~$(n+1)(n+2)=2\binom{n+2}{2}$.
Hence~$\ML(\sfT)$ divides~$H(\sfT)!$ for trees with~$\bd(\sfT)\le1$.
However, for~$\sfB_{2,2}$, the factor~$13$ of~$\ML(\sfB_{2,2})=52$ does not divide~$H(\sfB_{2,2})!=9!$.
Thus no formula that writes~$\ML(\sfT)$ as~$H(\sfT)!$ divided by a product of positive integer hook lengths holds for all trees of branching depth at most~$2$.

\section{Open problems}\label{sec:open}

Counting linear extensions is~$\#\mathrm{P}$-complete~\cite{BrightwellWinkler91}, even for posets with no chain of three elements~\cite[Theorem~1]{DittmerPak20}.

\begin{question}\label{q:complexity}
    Given a rooted tree~$\sfT$, can~$\ML(\sfT)$ be computed in polynomial time, or is this computation~$\#\mathrm{P}$-hard?
\end{question}

Explicit formulas for~$\ML(\sfB_{d,k})$ are known for~$k\le1$ (Proposition~\ref{prop:paths}) and for~$d\le2$ (Corollaries~\ref{cor:depth-two} and~\ref{cor:medium-brooms}); the case~$d\ge3$ and~$k\ge2$ remains open.

\begin{question}\label{q:broom-hook}
    Does Corollary~\ref{cor:medium-brooms} extend to a formula for~$\ML(\sfB_{d,k})$ whose number of summations is bounded independently of~$d$?
\end{question}

One possible approach is to seek a determinantal hook formula of the type obtained for mobile posets by Garver, Grosser, Matherne, and Morales~\cite{GGMM21}; Park~\cite{Park22} gave related Naruse-type hook formulas.

Our computations show that, for~$1\le n\le12$, every maximizer of~$\ML$ among rooted trees on~$n$ vertices is a broom.
The minimum is attained only by the path for~$n\le5$, by the spider~$\sfSp(3,1,1)$ for~$n=6$, and by the star for~$7\le n\le12$.

\begin{question}\label{q:extremal-ML}
    Which rooted trees~$\sfT$ on~$n$ vertices maximize or minimize~$\ML(\sfT)$?
    Does a broom attain the maximum for every~$n$?
\end{question}

\section*{Acknowledgments}
This project grew out of the Laboratory of Geometry at Michigan (LoGM) program at the University of Michigan.
It began in January 2025, and the core results were presented in a poster session organized by LoGM in April 2025.

The authors thank Andrew Sack and Katie Waddle for their mentorship and insightful suggestions.
They also thank Alejandro Bravo-Doddoli for helpful feedback.

All code was first written by hand.
Later, GPT Sol, GPT Astra, Claude Opus, and Claude Fable assisted with code generation, brute-force checking, and editing.

\bibliographystyle{amsplain}

\bibliography{references}

\end{document}